\documentclass{amsart}
\usepackage{geometry,amssymb,color,hyperref}
\usepackage{graphicx} 
\usepackage{mathtools}

\theoremstyle{plain}

\theoremstyle{definition}
\newtheorem{definition}{Definition}

\theoremstyle{remark}

\newcommand{\IZ}{\mathbb Z}

\title{There exist Steiner systems $S(2,7,505)$ and $S(2,7,589)$}
\author{Ivan Hetman}
\date{April 2026}

\begin{document}
    \begin{abstract} In this note two Steiner systems $S(2,7,505)$ and two Steiner systems $S(2,7,589)$ are presented.
    This resolves two of $21$ undecided cases for block designs with block length $7$, mentioned in Handbook of Combinatorial Designs.
    \end{abstract}

    \maketitle

    Given three integers $t$, $k$, $v$ such that $2 \le t < k < v$, a Steiner system $S(t, k, v)$ is a set $V$ of cardinality $v$ together with a family $\mathcal B$ of $k$-subsets of $V$ (called \textit{blocks}) such that every $t$-element subset of $V$ is contained in exactly one block.
    Handbook of Combinatorial Designs~\cite[II.3.5-6,8-9]{HoCD} lists $21$ numbers $v$ for which the existence of Steiner systems $S(2,7,v)$ is unknown.
    In this note we present two Steiner systems $S(2,7,505)$ and two Steiner systems $S(2,7,589)$ thus resolving two of the above undecided cases.

    Difference families are convenient and compact way to present Steiner systems.
    By a \emph{difference family} on a group $G$ of cardinality $v$ we understand a family $\mathcal B$ of $k$-element subsets (called \emph{base blocks}) of $G$ such that the pair $(G, \{Bg:B\in\mathcal B,g\in G\})$ is a Steiner system $S(2,k,v)$.
    For a base block $B$ denote by $\Delta B=\{ac^{-1}:a,c\in B,\;a\not= c\}$ its set of differences.
    It is well known that for coprime numbers $v$ and $k$, a family $\mathcal B$ of $k$-element subsets of a group of cardinality $v$ is a difference family if and only if the family of differences $\{\Delta B:B\in\mathcal B\}$ consists of pairwise disjoint sets of the cardinality $k(k-1)$.
    Using this property, it is easy to verify that two families presented below indeed are difference families on the cyclic group $\IZ_{505}=\{0,1,\dots,504\}$, producing the Steiner systems $S(2,7,505)$.
    Recognizing isomorphic designs and calculating automorphism groups are more complicated and require a specialized software like \texttt{nauty}~\cite{nauty} or calculating fingerprints introduced in~\cite{Het}.
    Difference families listed below are preceded with the order of the design automorphism group and the fingerprint witnessing that these designs are not isomorphic.

    \begin{enumerate}
    \item $2525\;\{1=15150, 2=424200, 3=14710650, 4=142803900, 5=475800900\} \newline \{\{0, 1, 3, 7, 47, 133, 284\}, \{0, 5, 70, 100, 173, 185, 476\}, \{0, 8, 82, 199, 248, 391, 474\}, \newline \{0, 9, 262, 298, 370, 386, 439\}, \{0, 10, 137, 156, 213, 233, 246\}, \{0, 11, 182, 250, 274, 414, 462\}, \newline \{0, 14, 78, 172, 366, 421, 477\}, \{0, 15, 136, 157, 174, 322, 344\}, \{0, 18, 181, 304, 330, 371, 442\}, \newline \{0, 23, 58, 120, 227, 287, 374\}, \{0, 25, 105, 150, 209, 260, 310\}, \{0, 27, 145, 206, 238, 416, 453\}\}$
    \smallskip
    \item $2525\;\{2=444400, 3=13483500, 4=139198200, 5=480628700\} \newline \{\{0, 1, 3, 7, 119, 242, 341\}, \{0, 5, 51, 63, 95, 254, 287\}, \{0, 8, 261, 297, 369, 388, 417\}, \newline \{0, 9, 26, 85, 357, 428, 490\}, \{0, 10, 170, 310, 391, 455, 492\}, \{0, 11, 93, 113, 279, 384, 422\}, \newline \{0, 14, 161, 290, 368, 421, 477\}, \{0, 16, 89, 158, 273, 327, 453\}, \{0, 18, 45, 79, 217, 304, 371\}, \newline \{0, 21, 43, 149, 206, 317, 456\}, \{0, 25, 135, 175, 328, 375, 450\}, \{0, 30, 227, 258, 324, 431, 470\}\}$
    \end{enumerate}

    It is possible to obtain more designs using mirroring technique.

    \begin{definition}
        A \emph{mirror} of a family $\mathcal B=\{B_1,\dots,B_n\}$ of subsets of a group is a family obtained by replacing some sets with their inverses.
        More formally, for any of $2^n$ functions $\mu:\{1,\dots,n\}\to \{1, -1\}$, the family $\mathcal B^\mu = \{B_1^{\mu(1)},\dots,B_n^{\mu(n)}\}$ is a mirror of $\mathcal B$.
    \end{definition}

    Computer calculations show that by mirroring from two presented above difference families it is possible to obtain $832$ non-isomorphic designs: $16$ with automorphism groups isomorphic to $\IZ_{505}\rtimes \IZ_5$, and $816$ with automorphism groups isomorphic to $\IZ_{505}$.

    \medskip

    Two examples below are difference families on the cyclic group $\IZ_{589}=\{0,1,\dots,588\}$.

    \begin{enumerate}
        \item $3534\;\{1=35340, 2=395808, 3=13199490, 4=190277628, 5=803917854\} \newline \{\{0, 1, 4, 258, 357, 455, 572\}, \{0, 2, 7, 223, 242, 401, 520\}, \{0, 6, 73, 102, 112, 393, 534\}, \newline \{0, 8, 251, 341, 426, 527, 555\}, \{0, 9, 103, 118, 316, 330, 523\}, \{0, 11, 95, 160, 187, 239, 369\}, \newline \{0, 12, 237, 267, 318, 390, 501\}, \{0, 13, 57, 120, 230, 509, 557\}, \{0, 16, 124, 177, 263, 443, 558\}, \newline \{0, 20, 136, 174, 439, 476, 553\}, \{0, 22, 68, 127, 272, 312, 463\}, \{0, 23, 58, 224, 284, 367, 525\}, \newline \{0, 24, 78, 121, 203, 421, 563\}, \{0, 25, 74, 165, 206, 425, 458\}\}$
        \item $3534\;\{2=296856, 3=14641362, 4=192398028, 5=800489874\} \newline \{\{0, 1, 95, 136, 269, 352, 580\}, \{0, 2, 26, 114, 317, 355, 515\}, \{0, 3, 89, 118, 256, 264, 482\}, \newline \{0, 4, 23, 56, 91, 399, 436\}, \{0, 5, 11, 75, 421, 457, 472\}, \{0, 7, 84, 254, 301, 323, 394\}, \newline \{0, 12, 166, 212, 275, 462, 534\}, \{0, 13, 62, 248, 420, 545, 572\}, \{0, 14, 222, 384, 426, 504, 549\}, \newline \{0, 16, 48, 108, 249, 408, 539\}, \{0, 18, 39, 119, 144, 396, 478\}, \{0, 20, 116, 240, 271, 305, 418\}, \newline \{0, 28, 81, 184, 242, 390, 487\}, \{0, 43, 164, 223, 313, 374, 447\}\}$
    \end{enumerate}

    Computer calculations show that by mirroring from two presented above difference families it is possible to obtain $2928$ non-isomorphic $S(2,7,589)$ designs: $16$ with automorphism groups or order $3934$, $24$ with automorphism groups of order $1767$, $336$ with automorphism groups of order $1178$, and $2552$ with automorphism groups of order $589$.

    \medskip

    In~\cite{CostaFengWang} single example of $S(2,8,624)$ was firstly found.
    In current paper all $10$ non-isomorphic designs generated by same group are presented.

    A slight modification of the constructions above produces $1$-\emph{rotational} designs (see~\cite[VI.6.74]{HoCD}) $S(2,8,624)$.
    In this case to the group $\IZ_{623}=\{0,1,\dots,622\}$ we add one more point $\infty$ which is fixed under the natural action of the group $\IZ_{623}$ on $\IZ_{623}\cup\{\infty\}$.
    For a subset $B$ of $\IZ_{623}\cup\{\infty\}$, $\Delta B=\{ac^{-1}:a,c\in B\cap \IZ_{623},\;a\not= c\}$.
    It is well known that a family $\mathcal B$ of $k$-element subsets of $\IZ_{v-1}\cup\{\infty\}$ is a difference family if and only if $\mathcal B = \{H\cup\{\infty\},B_1,\dots,B_n\}$, where $H$ is a subgroup of $G$, $\{\Delta B:B\in\mathcal B\}$ consists of pairwise disjoint sets, and $|\Delta B_i|=k(k-1)$ for all $i=1,\dots,n$.
    Each base block $B$ below produces the difference family $\mathcal B$ consisting of $H\cup\{\infty\}$, where $H = \{0, 89, 178, 267, 356, 445, 534\}$, and $11$ base blocks $B\cdot 8^i$ for $i=0,\dots,10$.
    The base block $H\cup\{\infty\}$ produces $89$ blocks of the design, and the other $11$ base blocks $B\cdot 8^i$ produce $623$ blocks of the design.
    The design isomorphic to one from~\cite{CostaFengWang} is third in the list.

    \begin{enumerate}
        \item $6853\;\{2=68530, 3=3673208, 4=55879362, 5=401119796, 6=976086496\} \newline B = \{0, 1, 3, 41, 216, 444, 462, 589\}$
        \item $6853\;\{2=68530, 3=3344264, 4=56290542, 5=393883028, 6=983241028\} \newline B =\{0, 1, 3, 118, 304, 350, 398, 435\}$
        \item $6853\;\{2=68530, 3=3508736, 4=56537250, 5=400626380, 6=976086496\} \newline B =\{0, 1, 3, 189, 286, 304, 568, 580\}$
        \item $6853\;\{2=68530, 3=3618384, 4=59826690, 5=393389612, 6=979924176\} \newline B =\{0, 1, 4, 11, 272, 343, 370, 519\}$
        \item $6853\;\{3=3453912, 4=53042220, 5=396213048, 6=984118212\} \newline B =\{0, 1, 4, 50, 384, 483, 571, 587\}$
        \item $6853\;\{2=68530, 3=4166624, 4=54070170, 5=396021164, 6=982500904\} \newline B =\{0, 1, 4, 197, 280, 335, 354, 601\}$
        \item $6853\;\{2=137060, 3=3728032, 4=56002716, 5=398954248, 6=978005336\} \newline B =\{0, 1, 4, 340, 434, 443, 471, 505\}$
        \item $6853\;\{2=137060, 3=3728032, 4=59621100, 5=397967416, 6=975373784\} \newline B =\{0, 1, 5, 35, 61, 177, 345, 414\}$
        \item $6853\;\{3=4440744, 4=57935262, 5=395472924, 6=978978462\} \newline B =\{0, 1, 7, 22, 62, 241, 276, 307\}$
        \item $6853\;\{3=4111800, 4=56783958, 5=398104476, 6=977827158\} \newline B =\{0, 1, 7, 207, 426, 463, 501, 531\}$
    \end{enumerate}

    Mirroring is also applicable to $1$-rotational difference families.
    By mirroring $10$ difference families above, it is possible to obtain $940$ non-isomorphic designs, $10$ with automorphism groups isomorphic to $\IZ_{623}\rtimes\IZ_{11}$, and $930$ designs with automorphism groups isomorphic to $\IZ_{623}$.

    \section*{Acknowledgement}

    I would like to express my sincere thanks to Taras Banakh and Alex Ravsky for valuable comments and improving design presentations.

\end{document}